\documentclass[12pt]{article}
\usepackage[numbers]{natbib}
\usepackage{graphicx}
\usepackage{color, multirow,amsfonts,amsmath}
\usepackage{url}
\usepackage{rotating}
\title{IAPO estimators in Exponentiated Fr$\acute{e}$chet case\\ (Work in progress)
\bigskip
\\Pavlina K. Jordanova
\\{\small{\it Faculty of Mathematics and Informatics, Konstantin Preslavsky University of Shumen, \\115 "Universitetska" str., 9712 Shumen, Bulgaria. \\Corresponding author:  pavlina\_kj@abv.bg.}}
\\Evelina Veleva
\\{\small{\it Department of Applied mathematics and Statistics, Angel Kanchev  University of Ruse, Bulgaria.}}
}

\begin{document}
\date{}

\maketitle

\begin{abstract}
In 2017 Jordanova and co-authors consider probabilities for $p$-outside values, and later on, they use them in order to construct distribution sensitive IPO estimators. These works do not take into account the asymmetry of the distribution. This shortcoming was recently overcome and the corresponding probabilities for asymmetric $p$-outside values, together with the so-called IAPO estimators, were defined. Here we apply these results to Exponentiated-Fr$\acute{e}$chet distribution, introduced in 2003 by Nadarajah and Kotz. The abbreviation "IAPO" comes from "Inverse Probabilities for Asymmetric P-Outside Values". These estimators use as an auxiliary characteristic the empirical asymmetric $p$-fences. In this way, the system relating the estimated parameters and the asymmetric probabilities for $p$-outside values has an easier solution. The comparison with our previous study about the corresponding IPO and IPO-NM estimators shows that IAPO estimators give better results for the index of regular variation of the right tail of the cumulative distribution function. A simulation study depicts their rates of convergence, and finishes this work. 
\end{abstract}

\section{INTRODUCTION}

Estimation of the tail parameters of heavy tailed distribution is useful for prediction of the quantiles outside the range of the data, and the expected return period. These could be, for example, Values at Risk (see e.g. Marinelli  et al. \cite{MarinelliCarlo2007}), future looses of some insurance company, damages caused by the next flood, or earthquake, or just the highest expected temperatures in some fixed region, and for some next period of time.

The Exponentiated Fr$\acute{e}$chet distribution is defined in Nadarajah and Kotz \cite{nadarajah2003exponentiated,nadarajah2006exponentiated}. It belongs to the set of heavy-tailed distributions. Its advantage is that its tail behaviour depends on two  parameters $\alpha > 0$ and $\lambda > 0$, and therefore, allows a better fit to the data. Here, analogously to \cite{jordanova2021ipo}, we define distribution sensitive estimator of the index of regular variation, which is appropriate not only for large but also for relatively small, and moderate samples. First, we compute probabilities for asymmetric $p$-outside values. Due to the Law of large numbers, they can be used for diagnostic of this distribution given a sample of independent observations. Then, we construct estimators, which are an alternative of the Maximum Likelihood estimator suggested in \cite{nadarajah2006exponentiated}. 

Consider a random variable (r.v.) $X$, with cumulative distribution function (c.d.f.) $F_X(x) = \mathbb{P}(X \leq x)$, $x \in \mathbb{R}$,  generalized left-continuous inverse (quantile function) $F_X^\leftarrow(y) = \inf\{x \in \mathbb{R}: F_X(x) \geq y\}$, $y \in (0, 1]$,
and  $p \in (0; 0.5]$. We assume that $F_X^\leftarrow(0) = \sup\{x \in \mathbb{R}: F_X(x) = 0\}$, and $\sup \emptyset = -\infty$.  The probabilities for asymmetric right $p$-outside values are defined in Jordanova \cite{MyPK2025} via the equality
\begin{equation}\label{PaRP}
p_{A,R,p}(X): = p_{A,R,p}(F_X):  = \mathbb{P}(X > R^A(F_X, p)), 
\end{equation}
where 
\begin{equation}\label{aRF}
R^A(F_X, p) = R^A(X,p) = F_X^\leftarrow(1 - p) + 2\frac{1-p}{p}(F_X^\leftarrow(1 - p) - F_X^\leftarrow(0.5)) 
\end{equation}
$$= \frac{2-p}{p}F_X^\leftarrow(1 - p) - 2\frac{1-p}{p}F_X^\leftarrow(0.5)= F_X^\leftarrow(0.5) + \frac{2-p}{p}(F_X^\leftarrow(1 - p) - F_X^\leftarrow(0.5))$$
is the theoretical asymmetric right $p$-fence. In this way, if the observed distribution is symmetric, $R^A(F_X, p)$ coincides with $R(F_X, p)$, defined in  \cite{Jordanova2020Monograph,JordanovaStehlikIPOestimation}, and $p_{A,R,p}(X)$ coincides with $p_{R,p}(X)$, introduced in the same works.
Analogously, probabilities for asymmetric left $p$-outside values are 
\begin{equation}\label{PaLP}
p_{A,L,p}(X): = p_{A,L,p}(F_X):  = \mathbb{P}(X < L^A(F_X, p)), 
\end{equation}
and the theoretical asymmetric left $p$-fence is
\begin{equation}\label{aLF}
L^A(F_X, p) = L^A(X,p) = F_X^\leftarrow(p) - 2\frac{1-p}{p}(F_X^\leftarrow(0.5) - F_X^\leftarrow(p)) 
\end{equation}
$$= \frac{2-p}{p}F_X^\leftarrow(p) - 2\frac{1-p}{p}F_X^\leftarrow(0.5) = F_X^\leftarrow(0.5) - \frac{2-p}{p}(F_X^\leftarrow(0.5) - F_X^\leftarrow(p)).$$
In the next section we will compute these characteristics of the Exponentiated Fr$\acute{e}$chet $l$-type with parameters $\alpha > 0$ and $\lambda > 0$. 

The figures in this work were plotted by using R software (2025) \cite{R}.

\section{PROBABILITIES FOR ASYMMETRIC $p$-OUTSIDE\\ VALUES}

Let us assume that a r.v. $X$ belongs to the Exponentiated Fr$\acute{e}$chet $l$-type with parameters $\mu \in \mathbb{R}$, $\sigma > 0$, $\alpha > 0$ and $\lambda > 0$.  According to the definition of Nadarajah and Kotz \cite{nadarajah2006exponentiated}, the last means that $X$ has a c.d.f. 
\begin{equation}\label{EFD}
F_X(x) =  1- \left[1-\exp\left\{-\left(\left(\frac{x - \mu}{\sigma}\right)^{-\lambda}\right)\right\}\right]^\alpha, \quad x > \mu,
\end{equation}
and $F_X(x) = 0$, otherwise. Briefly,  we will denote this fact by $X \in Exp-Fr(\alpha, \lambda; \mu, \sigma)$. As far as all probabilities for asymmetric right $p$-outside values are invariant with respect to shifting, and product with a positive scale parameter, without lost of generality, we assume that $\mu = 0$, and $\sigma = 1$.  Nadarajah and Kotz \cite{nadarajah2006exponentiated} show that this c.d.f. has a regularly varying right tail with parameter $-\alpha\lambda$. 
They use numerical methods in order to solve the Maximum Likelihood system of equations, and to estimate the parameters $\alpha$ and $\lambda$.  Here we provide an alternative of their approach, and construct estimators of $\alpha$ and $\lambda$, which are location and scale invariant.

{\it Note:} It is easy to see that for $\alpha = 1$ formula (\ref{EFD}) defines exactly the Fr$\acute{e}$chet distribution. Therefore, our results complement those considered e.g. in Gumbel \cite{gumbel1965quick} or Ramos et al. \cite{ramos2020frechet}.

The Exponentiated Fr$\acute{e}$chet quantile function is $F_X^\leftarrow(p) = \{-\log[1-(1-p)^{\frac{1}{\alpha}}]\}^{-\frac{1}{\lambda}}$ (see e.g. Rao et al. \cite{rao2016new}, or Jordanova \cite{Jordanova2020Monograph}). Therefore, formulae (\ref{aLF}) and (\ref{aRF}) lead us to 
\begin{eqnarray}
 \label{EFaLF} L^A(X,p)  & = & \frac{2-p}{p}\left\{-\log\left[1-(1-p)^{\frac{1}{\alpha}}\right]\right\}^{-\frac{1}{\lambda}} - 2\frac{1-p}{p}\left[-\log\left(1-2^{-\frac{1}{\alpha}}\right)\right]^{-\frac{1}{\lambda}}\\
\label{EFaRF}  R^A(X,p)  & = & \frac{2-p}{p}\left[-\log\left(1-p^{\frac{1}{\alpha}}\right)\right]^{-\frac{1}{\lambda}} - 2\frac{1-p}{p}\left[-\log\left(1-2^{-\frac{1}{\alpha}}\right)\right]^{-\frac{1}{\lambda}}.
  \end{eqnarray}
  
In this way we proved the following result.

{\bf Lemma 1.} If $X \in Exp-Fr(\alpha, \lambda; 0, 1)$, then for $p \in (0; 0.5]$,
\begin{itemize}
  \item[i)] $L^A(X, p)$ is presented in (\ref{EFaLF});
  \item[ii)] $R^A(X, p)$ is presented in (\ref{EFaRF}).
\end{itemize}
  
In order to compute probabilities for asymmetric left $p$-outside values, let us first see that $L^A(X,p) \geq 0$, if and only if, 
$$\frac{2-p}{2p(1-p)}\left\{-\log\left[1-(1-p)^{\frac{1}{\alpha}}\right]\right\}^{-\frac{1}{\lambda}} \geq \left[-\log\left(1-2^{-\frac{1}{\alpha}}\right)\right]^{-\frac{1}{\lambda}}$$
$$\iff \left(\frac{2p(1-p)}{2-p}\right)^\lambda\log\left[1-(1-p)^{\frac{1}{\alpha}}\right] \geq \log\left(1-2^{-\frac{1}{\alpha}}\right).$$
Therefore, if we denote by $p_0 \in (0, 0.5]$ the solution of the equation $\left(\frac{2p(1-p)}{2-p}\right)^\lambda\log\left[1-(1-p)^{\frac{1}{\alpha}}\right] =
 \log\left[1-2^{-\frac{1}{\alpha}}\right]$, by formulae (\ref{EFaLF}) and (\ref{PaLP}) we obtain that, $p_{A,L,p}(X) = 0$, when $p \in (0, p_0]$, and 

$p_{A,L,p}(X) = $
\begin{equation}\label{paL}
1 - \left\{1 - \exp\left\{-\left(\frac{2-p}{p}\left\{-\log\left[1-(1-p)^{\frac{1}{\alpha}}\right]\right\}^{-\frac{1}{\lambda}} - 2\frac{1-p}{p}\left\{-\log\left[1-2^{-\frac{1}{\alpha}}\right]\right\}^{-\frac{1}{\lambda}}\right)^{-\lambda}\right\}\right\}^\alpha 
   \end{equation}  
   if $p \in (p_0, 0.5]$.

Now, let us compute probabilities for asymmetric right $p$-outside values. It is easy to see that for all  $p \in (0; 0.5]$, $\frac{2-p}{p} \geq 2\frac{1-p}{p} > 0$; for all $\alpha > 0$, $1-p^{\frac{1}{\alpha}}\in (0, 1]$, $1-2^{-\frac{1}{\alpha}} \in (0, 1]$,  and $1-p^{\frac{1}{\alpha}} \geq 1-2^{-\frac{1}{\alpha}}$; and for all $\lambda > 0$, $h(x):=\left\{-\log(x)\right\}^{-\frac{1}{\lambda}}$ is an increasing function in $(0, 1]$. Then, 
$$\frac{2-p}{p}\left[-\log\left(1-p^{\frac{1}{\alpha}}\right)\right]^{-\frac{1}{\lambda}} \geq 2\frac{1-p}{p}\left[-\log\left(1-2^{-\frac{1}{\alpha}}\right)\right]^{-\frac{1}{\lambda}}.$$
The last entails that $R^A(X,p) \geq 0$ on the whole parametric space. Therefore, by using (\ref{PaRP}) and (\ref{EFaRF}) we obtain that for any $\lambda > 0$, $\alpha > 0$ and $p \in (0; 0.5]$, 
 \begin{equation}\label{paR}
p_{A,R,p}(X) =  \left\{1 - \exp\left\{-\left(\frac{2-p}{p}\left[-\log\left(1-p^{\frac{1}{\alpha}}\right)\right]^{-\frac{1}{\lambda}} - 2\frac{1-p}{p}\left[-\log\left(1-2^{-\frac{1}{\alpha}}\right)\right]^{-\frac{1}{\lambda}}\right)^{-\lambda}\right\}\right\}^\alpha.
\end{equation}

The above considerations prove the following result.

{\bf Theorem 1.} If $X \in Exp-Fr(\alpha, \lambda; \mu, \sigma)$, then for all $p \in (0; 0.5]$,
\begin{itemize}
\item[i)] $p_{A,L,p}(X)$ is presented in (\ref{paL}),
\item[ii)] $p_{A,R,p}(X)$ is presented in (\ref{paR}).
\end{itemize}

\begin{figure}[h]
\begin{minipage}[t]{0.5\linewidth}
\centerline{\includegraphics[width=1\textwidth]{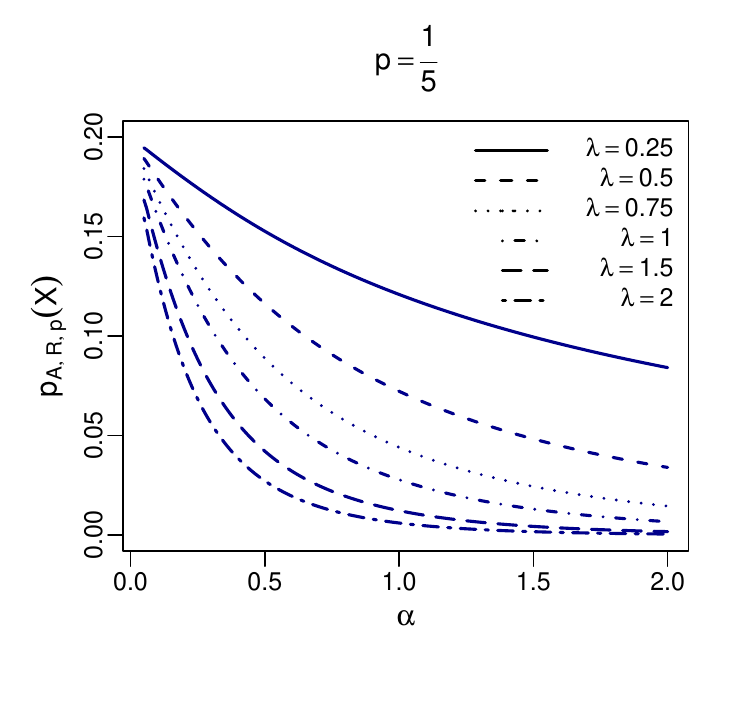}}
\end{minipage}
\begin{minipage}[t]{0.5\linewidth}
\centerline{\includegraphics[width=1\textwidth]{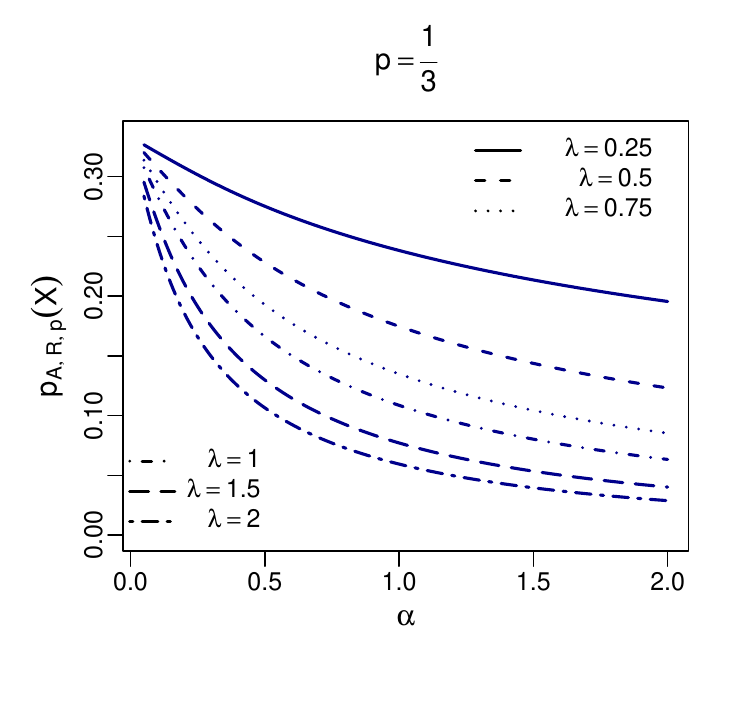}}
\end{minipage}
\caption{$p_{A,R,p}(X)$ in the Exponentiated Fr$\acute{e}$chet case for different values of $\alpha$ and $\lambda$. \label{fig:ExponentiatedFrechet1}}
\end{figure}

\begin{figure}[h]
\begin{minipage}[t]{0.5\linewidth}
\centerline{\includegraphics[width=1\textwidth]{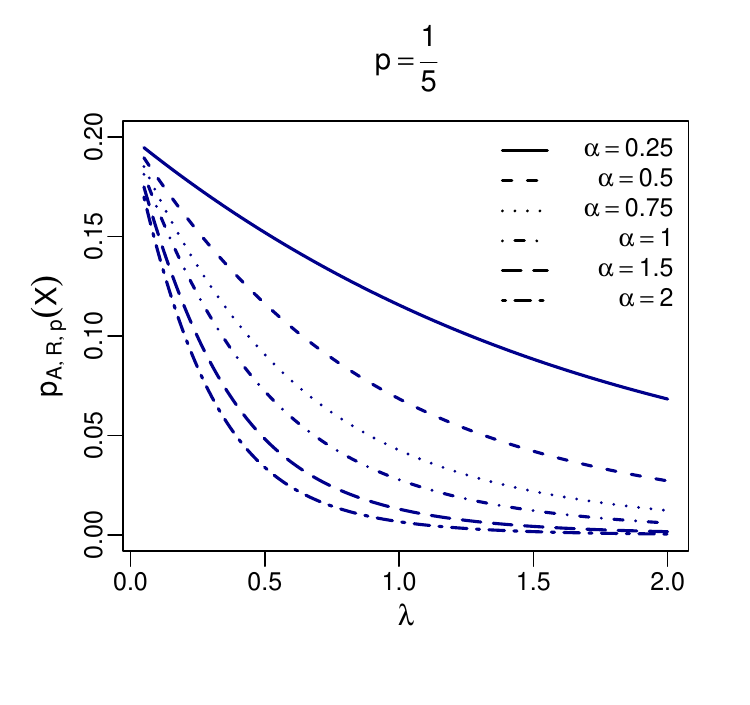}}
\end{minipage}
\begin{minipage}[t]{0.5\linewidth}
\centerline{\includegraphics[width=1\textwidth]{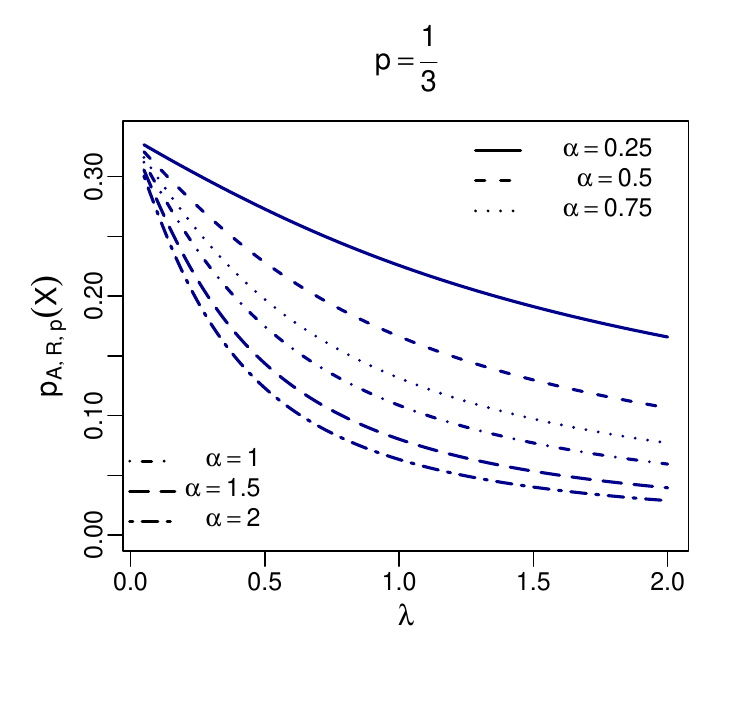}}
\end{minipage}
\caption{$p_{A,R,p}(X)$ in the Exponentiated Fr$\acute{e}$chet case for different values of $\alpha$ and $\lambda$. \label{fig:ExponentiatedFrechet2}}
\end{figure}

As we can see, the heaviness of the right tails depends on two parameters. Nadarajah and Kotz \cite{nadarajah2006exponentiated} showed that this dependence is via their product, and therefore, it is  symmetric with respect to $\alpha$ and $\lambda$. Figures \ref{fig:ExponentiatedFrechet1}-\ref{fig:ExponentiatedFrechet2} depict the dependencies of these probabilities on $\lambda$ and $\alpha$, for $p = \frac{1}{3}$, and $p = \frac{1}{5}$. When $\alpha$ decrease, the probabilities for asymmetric right $p$-outside values increase, together with the chance the considered r.v. to bring surprises. The dependence of $\lambda$ is analogous. Although, formula (\ref{paR}), is not symmetric with respect to $\lambda$ and $\alpha$, in these two figures we can observe similarities of the plots of $p_{A,R,p}(X)$  as functions of $\alpha > 0$, and as functions of $\lambda > 0$. 


\section{IAPO ESTIMATORS}
Consider  $p \in (0; 0.5]$, and independent identically distributed r.vs $X_1, X_2, ..., X_n$, $X \in Exp-Fr(\alpha, \lambda; \mu, \sigma)$. Their c.d.f. is defined in (\ref{EFD}).  In this section we interpret them as observations on $X$, and find new estimators of the parameters $\alpha$ and $\lambda$. 

Denote by  $X_{1:n} \leq X_{2:n} \leq \ldots X_{n:n}$  the corresponding increasing order statistics, and by 
\begin{equation}\label{ECF}
F_n(x):= \left\{\begin{array}{ccc}
                                    0 & , & x < X_{1:n} \\
                                    \frac{i}{n} & , & x \in [X_{i:n};X_{i+1:n}) \\
                                    1 & , & x \geq X_{n:n} 
                                  \end{array}
\right., \quad x \in \mathbb{R}.
\end{equation}
the empirical c.d.f. Then, $F_n\left(X_{i:n}\right) = \frac{i}{n}$, $i = 1, 2, \ldots, n$. There are many different definitions for the empirical $p$-quantiles. For $p \in (0; 1]$ we use the following one $F_n^\leftarrow(0) = X_{1:n}$, and for $p \in \left(\frac{i-1}{n}; \frac{i}{n}\right]$, i = 1, 2, \ldots, n,
\begin{equation}\label{quantiles}
F_n^\leftarrow (p): = \inf\{x \in \mathbb{R}, F_n(x) \geq p\} = X_{\lceil np \rceil :n} = X_{i:n}, 
\end{equation}
where $\lceil a \rceil$ means ceiling of $a$.  See for example \cite{serfling2009approximation,HyndmanandFan1996,Jordanova2020Monograph}. This function is left continuous and $F_n^\leftarrow \left(\frac{i}{n}\right) = X_{i:n}$, $i = 1, 2, \ldots, n$.
It coincides with 
\begin{equation}\label{empquant3}
F_n^\leftarrow(p) = \left\{\begin{array}{ccc}
                             X_{([np+1], n)} & , & np \not\in  \mathbb{N}, \\
                             X_{(np, n)} & , & np \in \mathbb{N},
                           \end{array}
\right. = \left\{\begin{array}{ccc}
                             X_{(\lceil np \rceil, n)} & , & np \not\in  \mathbb{N}, \\
                             X_{(np, n)} & , & np \in \mathbb{N},
                           \end{array}
\right.,
\end{equation}
 where $[a]$ is the integer part of $a$. The last definition is implemented by software R \cite{R} in function $quantile$, with parameter $Type = 1$.

Following Jordanova \cite{MyPK2025}, for $p \in (0; 0.5]$, let us denote by 
\begin{eqnarray*}
  R_n^A(p)  &=& F^\leftarrow_n(1 - p)  + 2\frac{1-p}{p}(F_n^\leftarrow(1 - p) - F_n^\leftarrow(0.5))  = \frac{2-p}{p}F_n^\leftarrow(1 - p) - 2\frac{1-p}{p}F_n^\leftarrow(0.5) \\
   &=& F_n^\leftarrow(0.5) + \frac{2-p}{p}(F_n^\leftarrow(1 - p) - F_n^\leftarrow(0.5)), \\
  L_n^A(p)  &=& F^\leftarrow_n(p) - 2\frac{1-p}{p}(F_n^\leftarrow(0.5) - F_n^\leftarrow(p)) 
= \frac{2-p}{p}F_n^\leftarrow(p) - 2\frac{1-p}{p}F_n^\leftarrow(0.5)\\
 &=& F_n^\leftarrow(0.5) - \frac{2-p}{p}(F_n^\leftarrow(0.5) - F_n^\leftarrow(p))
\end{eqnarray*}
the {\bf asymmetric empirical right} and {\bf asymmetric empirical left $p$-fences}.
Then, an observation $X$ is called {\bf sample (empirical) asymmetric right $p$-outside value} if $X > R_n^A(p)$, and {\bf sample (empirical) asymmetric left $p$-outside value} if $X < L_n^A(p).$ Let us denote by $n_{R}^A(p, n)$ the number of right asymmetric $p$-outside values in the sample, and by $\hat{p}_{R}^A (p, n) = \frac{n_{R}^A(p, n)}{n}$  their relative frequencies. Analogously to the algorithm introduced in Jordanova and Petkova (2018) \cite{MoniPoli2018}, for $p = 0.25$ and in Jordanova and Stehl{\'\i}k (2019) \cite{JordanovaStehlikIPOestimation} for general $p \in (0; 0.5]$ for IPO  estimators of the parameters $\alpha$ and $\lambda$, here we construct IAPO and IAPO-NM estimators. 
 For any fixed $0 < p_1 < p_2 < 0.5$ such that  $0 < \hat{p}_R^A(p_1, n) < \hat{p}_R^A(p_2, n)$, IAPO-NM estimators are  the solutions  $\hat{\hat{\alpha}}_{n}^A$, and $\hat{\hat{\lambda}}_{n}^A$ of the system  of equations
\begin{equation}\label{thesystem}
\left|\begin{array}{ccc}
          \hat{p}_R^A(p_1, n) & = & \left\{1 - \exp\left\{-\left(\frac{2-p_1}{p_1}\left[-\log\left(1-p_1^{\frac{1}{\alpha}}\right)\right]^{-\frac{1}{\lambda}} - 2\frac{1-p_1}{p_1}\left[-\log\left(1-2^{-\frac{1}{\alpha}}\right)\right]^{-\frac{1}{\lambda}}\right)^{-\lambda}\right\}\right\}^\alpha\\
           \hat{p}_R^A(p_2, n) & = &  \left\{1 - \exp\left\{-\left(\frac{2-p_2}{p_2}\left[-\log\left(1-p_2^{\frac{1}{\alpha}}\right)\right]^{-\frac{1}{\lambda}} - 2\frac{1-p_2}{p_2}\left[-\log\left(1-2^{-\frac{1}{\alpha}}\right)\right]^{-\frac{1}{\lambda}}\right)^{-\lambda}\right\}\right\}^\alpha
        \end{array}\right..
\end{equation}

Our experiments show that the solution of (\ref{thesystem}) is too sensitive to accuracy of the estimators $\hat{p}_R^A(p_1, n)$ and $\hat{p}_R^A(p_2, n)$. This disadvantage is overcome via the next estimators. The next step when work with IAPO-NM estimators is to estimate the shift parameter $\mu$, and the scale parameter $\sigma$ described in (\ref{EFD}). We can make this via the quantile matching procedure which is described e.g. in Klugman et al. (2012) \cite{klugman2012loss} or in Sgouropoulos et al. (2015) \cite{sgouropoulos2015matching}. In this way we have divided the estimation procedure of the four parameters $\alpha$, $\lambda$, $\mu$ and $\sigma$ in two system  of equations: one for $\alpha$ and $\lambda$, and after that one for $\mu$ and $\sigma$.

 Let us now consider IAPO estimators. The last abbreviation comes from the expression "Inverse Asymmetric Probabilities for P-Outside values". 
Now, we need some preliminary information about the parameters $\mu$ and $\sigma$ in (\ref{EFD}). Given $\mu$, $\sigma$, and $0 < \hat{p}_R^A(p_1, n) < \hat{p}_R^A(p_2, n)$, we obtain IAPO statistics $\hat{\alpha}_{n}^A$, and $\hat{\lambda}_{n}^A$ as the real positive solution, of the system of equations
\begin{equation}\label{System}
\left|\begin{array}{ccc}

  \hat{p}_R^A(p_1, n) &=& \left\{1 - \exp\left(-\left(\frac{R_n^A(p_1)-\mu}{\sigma}\right)^{-\lambda}\right)\right\}^{\alpha} \\
  \hat{p}_R^A(p_2, n) &=& \left\{1 - \exp\left(-\left(\frac{R_n^A(p_2)-\mu}{\sigma}\right)^{-\lambda}\right)\right\}^{\alpha}
\end{array}\right..
\end{equation}

$$\left|\begin{array}{ccc}

  \log\left(\hat{p}_R^A(p_1, n)\right) &=& \alpha\log\left\{1 - \exp\left(-\left(\frac{R_n^A(p_1)-\mu}{\sigma}\right)^{-\lambda}\right)\right\} \\
  \log\left(\hat{p}_R^A(p_2, n)\right) &=& \alpha\log\left\{1 - \exp\left(-\left(\frac{R_n^A(p_2)-\mu}{\sigma}\right)^{-\lambda}\right)\right\}
\end{array}\right..$$
Thus, 
$$c = \frac{\log\left\{1 - \exp\left(-\left(\frac{R_n^A(p_1)-\mu}{\sigma}\right)^{-\lambda}\right)\right\}}{\log\left\{1 - \exp\left(-\left(\frac{R_n^A(p_2)-\mu}{\sigma}\right)^{-\lambda}\right)\right\}}$$
where $c = \log_{\hat{p}_R^A(p_2, n)}\left(\hat{p}_R^A(p_1, n)\right)$, and
\begin{equation}\label{EqIAPOlambda}
\left\{1-\exp\left(-\left(\frac{R_n^A(p_2)-\mu}{\sigma}\right)^{-\lambda}\right)\right\}^c = 1-\exp\left\{-\left(\frac{R_n^A(p_1)-\mu}{\sigma}\right)^{-\lambda}\right\}.
\end{equation}
The numerical solution of the last equation with respect to $\lambda$ is exactly the IAPO estimator $\hat{\lambda}_{n}^A$. 
Now, the IAPO estimator for $\alpha$ is
\begin{equation}\label{EqIAPOalpha}
\hat{\alpha}_{n}^A = \frac{\log\left(\hat{p}_R^A(p_1, n)\right)}{\log\left\{1 - \exp\left(-\left(\frac{R_n^A(p_1)-\mu}{\sigma}\right)^{-\lambda}\right)\right\}}.
\end{equation}
 Let us remind that these estimators are distribution sensitive. Therefore, only if our assumption $X \in Exp-Fr(\alpha, \lambda; \mu, \sigma)$ is correct, then, the above expression would not depend substentially on the choise of $p_1$. 

 The corresponding estimator of the quantiles outside the range of the data is
$$\hat{F}^\leftarrow_{n}(p) = \mu + \sigma\left\{-\log[1 - (1 - p)^{1/\hat{\alpha}_{n}^A}]\right\}^{-1/\hat{\lambda}_{n}^A}.$$

\section{SIMULATION STUDY}

\begin{figure}
\begin{minipage}[t]{0.5\linewidth}
\includegraphics[scale=.55]{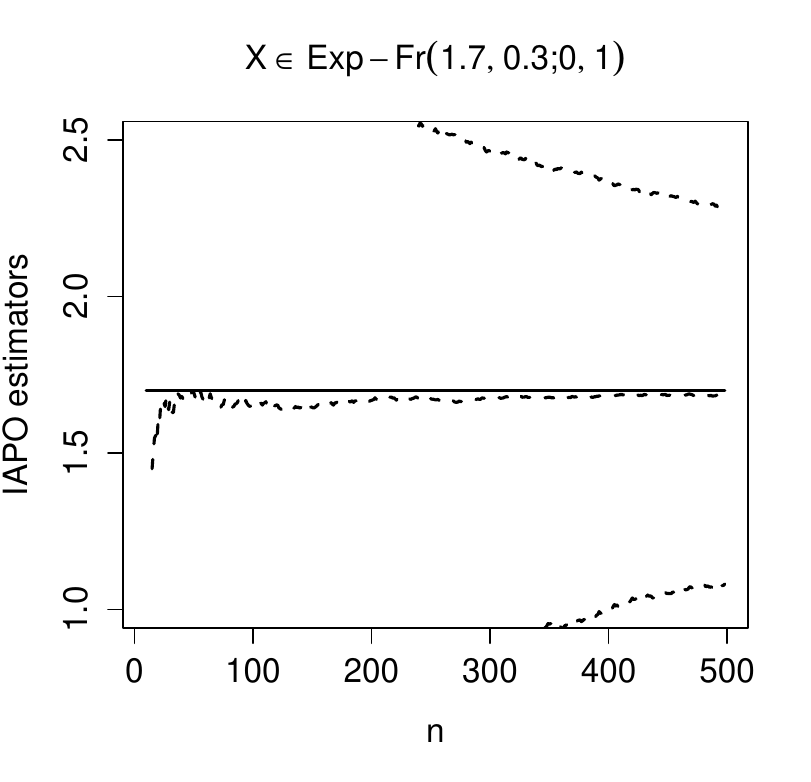}\vspace{-0.5cm}
\end{minipage}$ $
\begin{minipage}[t]{0.5\linewidth}
\includegraphics[scale=.55]{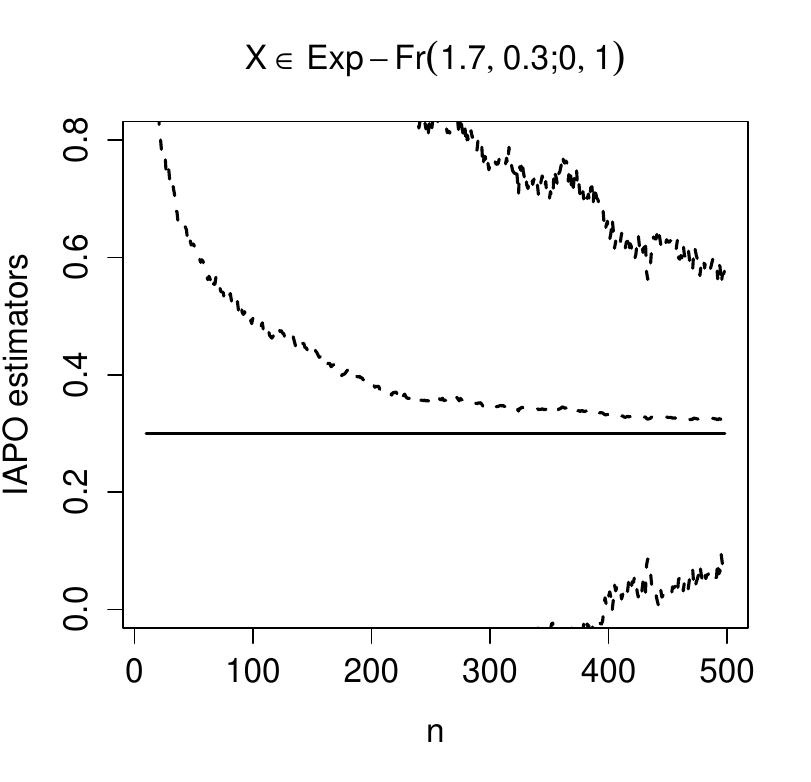}\vspace{-0.5cm}
\end{minipage}
\caption{\small IAPO estimators: Dependence of $\hat{\alpha}_{n}^A$(left) and $\hat{\lambda}_{n}^A$(right) on $n$ for $p_1 = \frac{1}{3}$ and $p_2 = \frac{1}{5}$.} \label{fig:SymStudyM10}
\end{figure}

\begin{center}
\begin{figure}
\begin{center}
\includegraphics[scale=.55]{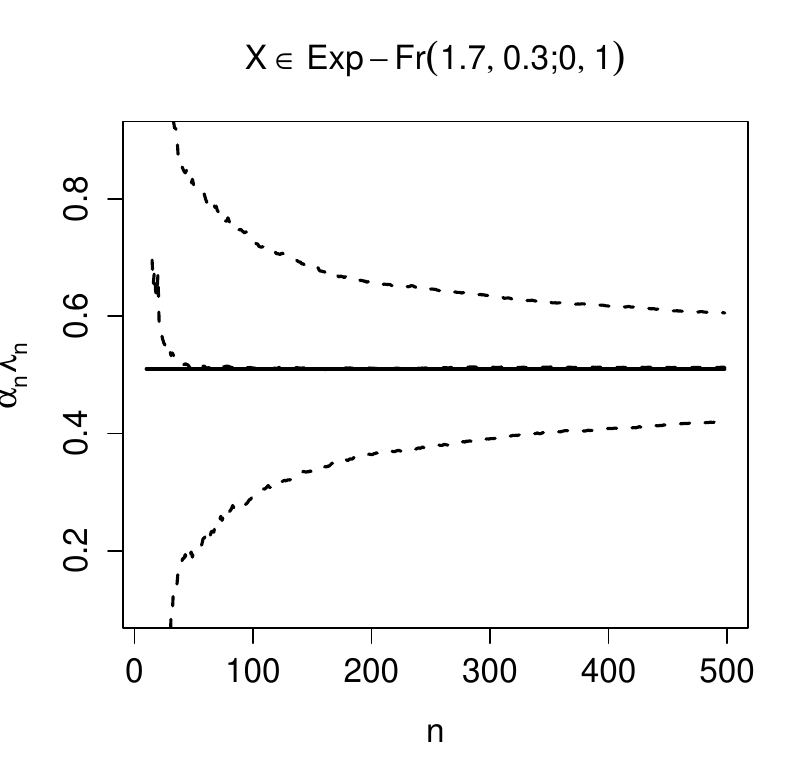}\vspace{-0.5cm}
\caption{\small Dependence of the product of IAPO estimators $\hat{\alpha}_{n}^A$ and $\hat{\lambda}_{n}^A$ on $n$ for $p_1 = \frac{1}{3}$ and $p_2 = \frac{1}{5}$.} \label{fig:SymStudyAL}
\end{center}
\end{figure}
\end{center}

In this section, first we simulate $1000$ samples of $500$
independent observations on a r.v. $X \in Exp-Fr(\alpha, \lambda;
0, 1)$, for fixed $\alpha > 0$ and $\lambda > 0$. Then, for $p_1 = \frac{1}{3}$, $p_2 = \frac{1}{5}$, for any
fixed $n = 30, 31, ..., 500$, and for any fixed sample we 
compute $\hat{\alpha}_{n}^A$, and $\hat{\lambda}_{n}^A$,
the corresponding averages and standard errors over the samples.
These choices of $p_1$, and $p_2$ were done on our experience with the data, and the results along to our simulation study. Theoretical evidence related to the best their values are still an open problem. When $p_1$, or the sample size, are too small, we have no empirical right asymmetric $p_1$-outside values, and therefore, IAPO estimators are not applicable. 
When $p_1$, and $p_2$ are too close, the difference between $p_{R, p_i}$, $i = 1, 2$ is too small, and we need to use a higher accuracy of the computations in order to obtain the solution of (\ref{System}). When $p_2 \approx 0.5$, then the asymmetric $p_2$ right fence will be close to the median, the probabilities for asymmetric $p_2$ outside values will be close to $0.5$, which is uninformative, does not depend on the type of the distribution, and leads to higher deviations of the estimators. Therefore, the choice of $p_1$ and $p_2$ as possible as different, and which belongs to the interval $(0; 0.5]$, but are not close neither to $0$ nor to $0.5$ seems to be the best one.
Let us now depict the quality of the IAPO estimators $\hat{\alpha}_{n}^A$ and $\hat{\lambda}_{n}$ in some particular cases.
\begin{figure}
\begin{minipage}[t]{0.5\linewidth}
\includegraphics[scale=.55]{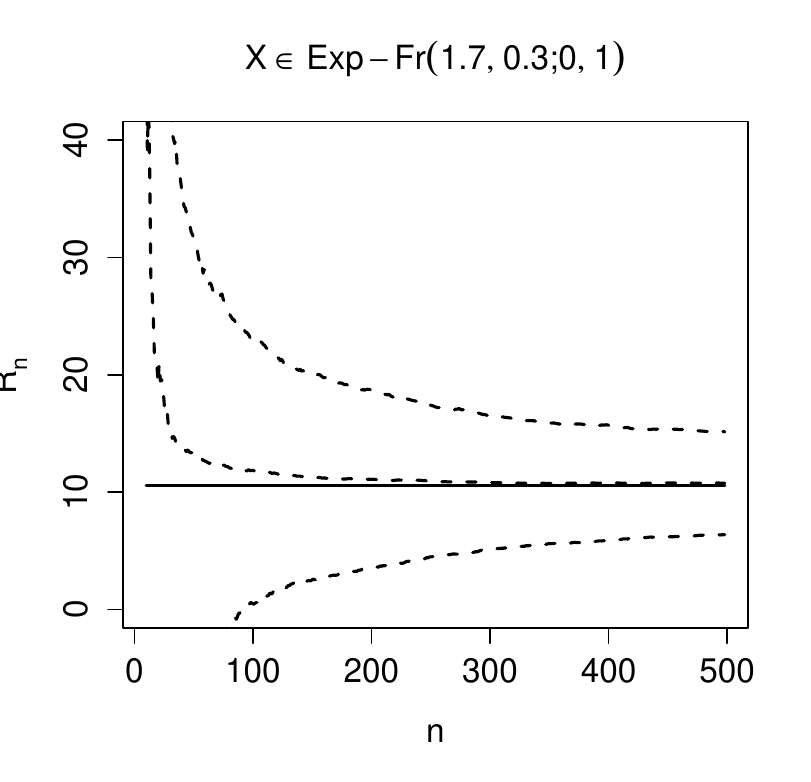}\vspace{-0.5cm}
\end{minipage}$ $
\begin{minipage}[t]{0.5\linewidth}
\includegraphics[scale=.55]{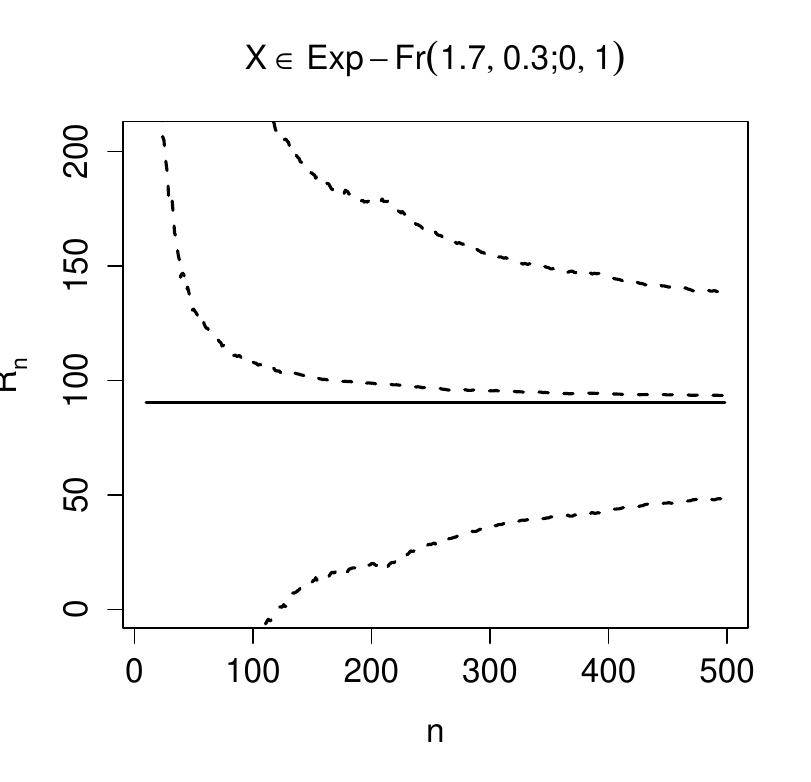}\vspace{-0.5cm}
\end{minipage}
\caption{\small Dependence of $R_n^A\left(\frac{1}{3}\right)$(left) and $R_n^A\left(\frac{1}{5}\right)$(right) on $n$.} \label{fig:SymStudyF}
\end{figure}
\begin{figure}
\begin{minipage}[t]{0.5\linewidth}
\includegraphics[scale=.55]{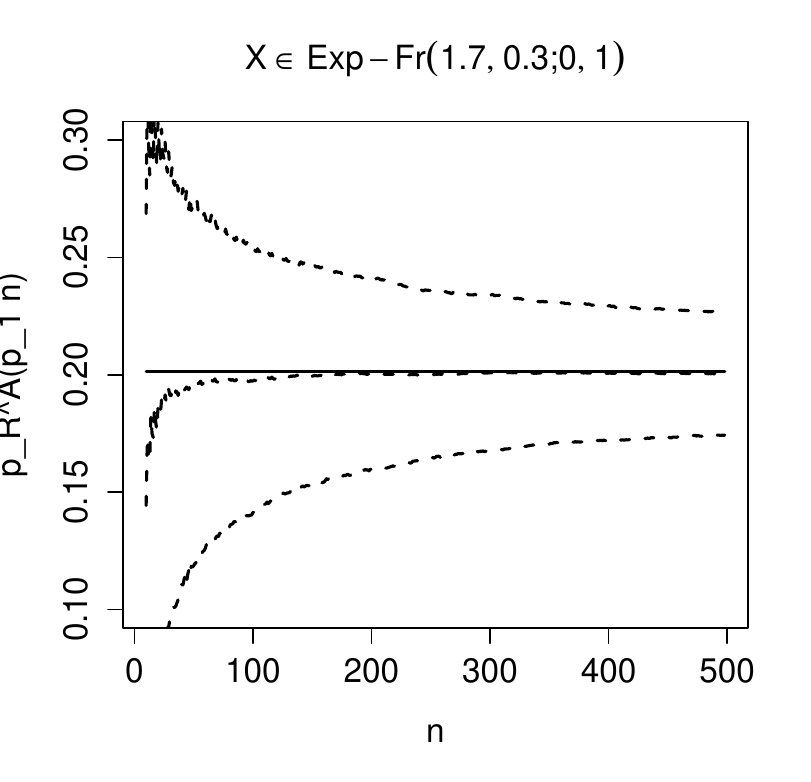}\vspace{-0.5cm}
\end{minipage}$ $
\begin{minipage}[t]{0.5\linewidth}
\includegraphics[scale=.55]{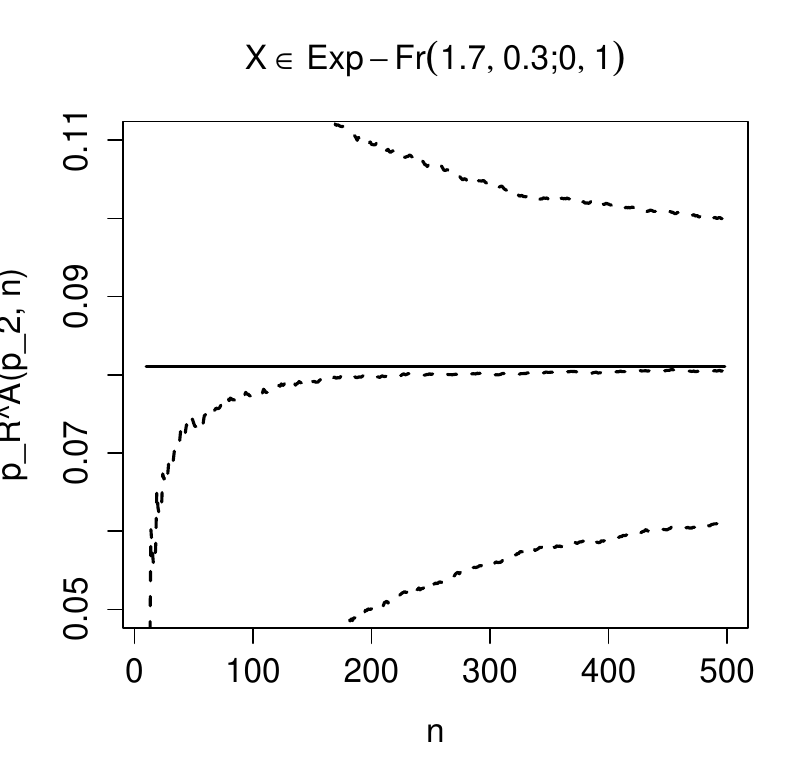}\vspace{-0.5cm}
\end{minipage}
\caption{\small Dependence of $\hat{p}_R^A\left(\frac{1}{3}, n\right)$(left) and $\hat{p}_R^A\left(\frac{1}{5}, n\right)$(right) on $n$.} \label{fig:SymStudyP}
\end{figure}
The values of $R_n^A\left(p\right)$, together with the corresponding asymptotically normal
$0.95$-confidence intervals are presented in Figure \ref{fig:SymStudyF}. The rate of convergence of $\hat{p}_R^A\left(p, n\right)$ estimators is depicted in Figure \ref{fig:SymStudyP}. The results for IAPO estimators of $\alpha$ and $\lambda$, together with the corresponding asymptotically normal
$0.95$-confidence intervals are plotted in Figure \ref{fig:SymStudyM10}. When speak about separate estimation of these two parameters we should mention that in some cases, for example when $\alpha = 0.5$ and $\lambda = 1.3$, these estimators give pure results. However, in any case the estimators of $\alpha \lambda$ have very nice properties because this algorithm is developed for estimation of the tail parameter, and for Exponentiated Fr$\acute{e}$chet distribution it is equal to $\alpha\lambda$. The rate of convergence of the estimator of the product is presented in particular in Figure \ref{fig:SymStudyAL}.  The solid lines show the values of the corresponding
estimated parameters. When compare these results with those in \cite{jordanova2021ipo} we observe that the results are very similar. This leads us to the hypothesis that IPO ansd IAPO estimators are independent on the best estimators of the asymmetry of the distribution. All these figures show that the bigger the sample size, the better  the quality of these estimators is. Similarly to IPO estimators these estimators are consistent. See Jordanova (2025) \cite{MyPK2025},
Jordanova and Stehlik (2019) \cite{JordanovaStehlikIPOestimation}, Jordanova (2020) \cite{Jordanova2020Monograph}, Jordanova and Nedzhibov \cite{jordanova2021ipo} among others.


\section{CONCLUSIVE REMARKS AND OPEN PROBLEMS}
The current work shows the usefulness of probabilities for asymmetric $p$-outside values for diagnostic, and estimation of the parameters which govern the tail behaviour of the Exponentiated Fr$\acute{e}$chet distribution. These universal characteristics, introduced in Jordanova \cite{MyPK2025} are invariant within linear probability type. They always exist. They outperform, for example, the role of the kurtosis. They refine the well-known concepts for heavy-tailed distribution, and have more general applications than the index of regular variation.

From statistical point of view, they allow us to construct new estimators of the parameters which govern the tail behaviour of the observed distribution. Although they are distribution sensitive, they give possibility to check if the observed distribution is possible to be Exponentiated Fr$\acute{e}$chet. When the observed sample is moderate or small they outperform many of the well-known non-parametric estimators of the index of regular variation, for example, Hill (1975), Pickands (1975) and Dekkers-Einmahl-de Haan Dekkers (1989) estimators. The formulae, and references about them could be found, for example in Embrechts et al. (2013) \cite{EKM}. The suggested algorithm gives better results for distributions with heavier tails.
These results are useful for practitioners who would like to estimate the quantiles outside the range of the data and to analyze extremal values.

\section{ACKNOWLEDGMENTS}
The first author was partially supported by the Project
RD-08-118/06.02.2025 from the Scientific Research Fund in
Konstantin Preslavsky University of Shumen, Bulgaria.


\nocite{*}
\bibliographystyle{plain}%

\begin{thebibliography}{99}
\bibitem{MarinelliCarlo2007} Marinelli, C.;  d\'Addona, S.;  Rachev, S., T., A comparison of some univariate models for {\MakeUppercase{V}}alue-at-{\MakeUppercase{R}}isk and expected shortfall, {\it{International Journal of Theoretical and Applied Finance}}, vol. 10(06), pp. 1043-1075, World Scientific Publishing, 2007.
\bibitem{nadarajah2003exponentiated} Nadarajah, S.; Kotz, S., The exponentiated {\MakeUppercase{F}}r{\'e}chet distribution, {\it{Interstat Electronic Journal}}, vol. 14, pp. 1-7, Citeseer, 2003.
\bibitem{nadarajah2006exponentiated} Nadarajah, S.; Kotz, S., The exponentiated type distributions, {\it{Acta Applicandae Mathematica}}, vol. 92, pp. 97--111, Springer, 2006.
\bibitem{jordanova2021ipo} Jordanova, P. K.; Nedzhibov, G. H., IPO and IPO-NM estimators in exponentiated Fr{\'e}chet case, {\it{AIP Conference Proceedings}},
  vol. 2333 (1), Article Number 150001, AIP Publishing, 2021.
\bibitem{MyPK2025} Jordanova, P.K., Probabilities for asymmetric $p$-outside values,
Mathematics and education on mathematics, {\it{ Proceedings of the Fifty-Fourth Spring Conference
of the Union of Bulgarian Mathematicians, Varna, March 31 – April 4, 2025}}, pp. 51-60, Union of Bulgarian Mathematicians, 2025.
  {\url{http://www.math.bas.bg/smb/2025_PK/tom_2025/pdf/051-060.pdf}}
\bibitem{gumbel1965quick} Gumbel, E.J., A quick estimation of the parameters in Fr{\'e}chet's distribution, {\it{Revue de l'Institut International de Statistique}},
  pp. 349--363, JSTOR, 1965.
\bibitem{ramos2020frechet} Ramos, P.L.; Louzada, Fr.; Ramos, Ed.; Dey, S., The Fr{\'e}chet distribution: Estimation and application-An overview,
  {\it{Journal of Statistics and Management Systems}}, vol. 23(3), pp. 549--578, Taylor \& Francis, 2020.
\bibitem{rao2016new} Rao, G.Sr.; Rosaiah, K.; Babu, M. Sr.; Sivakumar, D. Ch., New acceptance sampling plans based on percentiles for exponentiated Fr{\'e}chet distribution,
  Economic Quality Control, vol. 31(1), pp. 37--44, De Gruyter, 2016.
\bibitem{Jordanova2020Monograph} Jordanova, P., {\it{Probabilities for $p$-outside values and heavy tails}}, Shumen University Publishing House, 2020.
{\url{https://pavlinakj.wordpress.com/2017/09/29/pavlina-jordanova/}}
\bibitem{serfling2009approximation} Serfling, R. J., Approximation theorems of mathematical statistics, vol. 162, {\it{John Wiley \& Sons}}, 2009.
\bibitem{HyndmanandFan1996} Hyndman, R.J.; Fan, Y., Sample quantiles in statistical packages, {\it{The American Statistician}}, vol. 50(4), pp. 361--365, Taylor \& Francis, 1996.
\bibitem{R} R Development Core Team, R: A Language and Environment for Statistical Computing, {\it{R Foundation for Statistical Computing}}
  {\url{https://www.r-project.org/}},
\bibitem{JordanovaStehlikIPOestimation} Jordanova, P.; Stehl{\'\i}k, M., {\MakeUppercase{IPO}} estimation of the heaviness of the distribution beyond regularly varying tails,
  {\it{Stochastic Analysis and applications}}, vol. 38(1), pp.76-96, Taylor \& Francis, 2020.
 {\url{https://www.tandfonline.com/doi/abs/10.1080/07362994.2019.1647786?journalCode=lsaa20}}
\bibitem{MoniPoli2018} Jordanova, P.;  Petkova, M., Tails and probabilities for extreme outliers,
{\it{AIP Conference Prococeedings}}, vol. 2025, pp. 030002-1--030002-9, USA: AIP Publishing LLC, 2018.
\bibitem{klugman2012loss}Klugman, St. A.; Panjer, H. H.; Willmot, G. E.; Loss models: from data to decisions, vol. 715, {\it{John Wiley \& Sons}}, 2012.
\bibitem{sgouropoulos2015matching} Sgouropoulos, N.; Yao, Q.; Yastremiz, Cl., Matching a distribution by matching quantiles estimation, {\it{Journal of the American Statistical Association}}, vol. 110(510), pp. 742-759, Taylor \& Francis, 2015.
\bibitem{EKM} Embrechts, P.; Kl{\"u}ppelberg, Cl.; Mikosch, Th., Modelling extremal events: for insurance and finance, vol. 33, Springer Science \& Business Media, 2013.

\end{thebibliography}

\end{document}